\documentclass [12pt, draft] {amsart}
\usepackage {amssymb, amsmath, amscd}

\theoremstyle {plain}
\newtheorem {theorem} {Theorem}

\newtheorem {lemma} {Lemma}

\theoremstyle {definition}
\newtheorem {definition} {Definition}

\theoremstyle {remark}
\newtheorem* {remark} {Remark}

\newtheorem* {acknowledgments} {Acknowledgments}

\numberwithin {equation} {section}

\begin {document}
\title {Some endomorphisms of the hyperfinite $II_1$ factor}
\author {Hsiang-Ping Huang}
\address {Department of Mathematics, University of Utah, Salt Lake City, UT 84112}
\email {hphuang@math.utah.edu}
\thanks {Research partially supported by National Center for Theoretical
         Sciences, Mathematics Division, Taiwan}
\keywords { binary shifts, endomorphisms, irreducibility, subfactor}
\subjclass [2000] {46L37, 47B47}

\pagestyle{plain}

\begin {abstract}
For any finite dimensional $C^*$-algebra $A$ with any trace vector
$\vec s$ whose components are rational numbers, we give an endomorphism
$\Phi$ of the hyperfinite $II_1$ factor $R$  such that:
\[
  \forall \ k \ \in \ \mathbb {N}, \
  \Phi^k (R)' \cap R= \otimes^k A.
\]
The canonical trace $\tau$ on $R$ extends the trace vector $\vec s$
on $A$.

As a corollary,
we construct a one-parameter family 
of inclusions of hyperfinite $II_1$ factors
$N^{\lambda} \subset M^{\lambda}$
with trivial relative commutant $(N^{\lambda})' \cap M^{\lambda}= {\mathbb C}$
and with the Jones index 
\[
 [M^{\lambda}: N^{\lambda}]= \lambda^{-1} 
 \in (4, \infty) \cap {\mathbb Q}
\]
This partially solves the problem of finding all possible values
of indices of subfactors with trivial relative commutant
in the hyperfinite $II_1$ factor, by showing that
any rational number $\lambda^{-1} > 4$ can occur.   
\end {abstract}
\maketitle
\tableofcontents
\section {Introduction} \label {S: intro}

Subfactor theory \cite {vJ} is to
describe the position of a subfactor $N$ in an ambient factor $M$.
The standard invariant 
associated with Jones basic construction ${\mathcal G}_{N, M}$, 
\[
 \begin{matrix}
     &N' \cap N  \subset  & N'  \cap M  \subset  & N' \cap M_1  \subset 
      & N' \cap M_2  \subset &\cdots\\
     & \cup        &\cup        & \cup               & \cup                &\ \\
     & M' \cap N  \subset  & M' \cap M  \subset   &M' \cap M_1 \subset 
      & N' \cap M_2  \subset & \cdots
 \end{matrix}
\]
is a complete invariant in the amenable case \cite{sP94}.
To classify the standard invariant is 
the central topic 
ever since V.Jones founded the subfactor theory.

For a hyperfinite $II_1$ subfactor of finite Jones index, 
it is equipped with an extra structure: an endomorphism $\Phi$,
sending the ambient factor $M$ onto the subfactor $N$.
Therefore it is only natural to investigate the role of
the endomorphism.

A well-known example of endomorphisms is the canonical shift in a strongly 
amenable inclusion.
Another surprising example is  the binary shift \cite {rP88}
which gives rise to a counterexample that fails the tensor product formula
for entropy \cite{hN95}.
Via the Cuntz algebra, a lot of endomorphisms have been manufactured.

The Jones index $[M: \Phi (M)]$ is an outer-conjugacy invariant for 
endomorphisms. 
In the case of finite Jones index,
there is a distinguished outer-conjugacy invariant: the tower
of inclusions of finite dimensional $C^*$ algebras, $\{ A_k= \Phi^k
(M)' \cap M \}_{k=1}^{\infty}$.
The main part of the paper is to investigate the above 
invariant.

We prove that
for any finite dimensional $C^*$-algebra $A$ with any trace vector
$\vec s$ whose components are rational numbers, there exists an endomorphism
$\Phi$ of the hyperfinite $II_1$ factor $R$  such that:
\[
  \forall \ k \ \in \ \mathbb {N}, \
  \Phi^k (R)' \cap R= \otimes^k A.
\]
The canonical trace $\tau$ on $R$ extends the trace vector $\vec s$
on $A$.

Due to the idiopathic behavior of Powers' binary shift, 
our main result has a unexpected feedback to its origin:
the classification of hyperfinite $II_1$ subfactors.
In short, there is an analogy between Powers' binary shift 
and free product with amalgation.

As an application, we partially solve the problem of finding 
all possible values of indices of subfactors with trivial relative 
commutant in the hyperfinite $II_1$ factor, by showing that
any rational number $\lambda^{-1} > 4$ can occur.

\section {Preliminaries} \label {S: prelim}
Let $M$ be a $II_1$ factor with the canonical trace $\tau$. Denote
the set of unital *-endomorphisms of $M$ by $End (M, \tau)$. Then
$\Phi \in End (M, \tau)$ preserves the trace and $\Phi$ is
injective. $\Phi (M)$ is a subfactor of $M$. If there exists a
$\sigma \in Aut (M)$ with $\Phi_1 \cdot \sigma= \sigma \cdot \Phi_2$
for $\Phi_i \in End (M, \tau)$ $(i= 1, 2)$ then $\Phi_1$ and
$\Phi_2$ are said to be conjugate. If there exists a $\sigma \in Aut
(M)$ and a unitary $u \in M$ such that $Ad u \cdot \Phi_1 \cdot
\sigma= \sigma \cdot \Phi_2$, then $\Phi_1$ and $\Phi_2$ are outer
conjugate.

The Jones index $[M: \Phi (M)]$ is an outer-conjugacy invariant. We
consider only the finite index case unless otherwise stated. In such
case, there is a distinguished outer-conjugacy invariant: the tower
of inclusions of finite dimensional $C^*$ algebras, $\{ A_k= \Phi^k
(M)' \cap M \}_{k=1}^{\infty}$.

\begin {lemma}
$A_k = \Phi^k (M)' \cap M$ contains an subalgebra that is isomorphic
to $\otimes_{i= 1}^k A_1$, the $k$-th tensor power of $A_1$, where
$A_1 = \Phi (M)' \cap M$ as denoted.
\end {lemma}

\begin {proof}
We collect some facts here, for $1 \leq i < k$: \\
(1) $\Phi^{i} (A_1)$ is isomorphic to $A_1$,
since $\Phi^{i}$ is injective.\\
(2) $ A_1 \cap \Phi^{i} (A_1) \subseteq \Phi (M)'
\cap \Phi (M) = \mathbb{C}$.\\
(3) $[A_1, \Phi^{i} (A_1)] = 0$ by the definition of $A_1$.\\
(4) $\Phi^{i} (A_1) \subseteq A_{i+1} \subseteq A_k$.\\
\end {proof}

\begin {remark}
The dimension of the relative commutant ${\Phi^k (M)}' \cap M$ is
known to be bounded above by the Jones index $[M: \Phi(M)]^k$. The 
above lemma provides the lower bound for the growth estimate.
\end {remark}

A good example is the canonical shift \cite {dB} on the tower of
higher relative commutants for a strongly amenable inclusion of
$II_1$ factors of finite index. The ascending union of higher
relative commutants gives the hyperfinite $II_1$ factor, and the
canonical shift can be viewed as a $*$-endomorphism on the hyperfinite
factor. Lemma 1 is nothing but the commutation relations in S.Popa's
$\lambda$-lattice axioms \cite {sP95}.

A natural question arises with the above observation: for any finite
dimensional C*-algebra $A$, can we find a $II_1$ factor $M$ and a
$\Phi \in End (M, \tau)$ such that for all $k \in {\mathbb N}$,
\[
  \Phi^k (M)' \cap M \simeq \otimes_{i=1}^k A_i,
  {\text {\ where \ }} A_i \simeq A?
\]

The answer is positive and furthermore we can choose
$M$ to be the hyperfinite $II_1$ factor $R$.
We give the construction
in the next section. The main technical tool in
the construction is \cite {rP88} R.Powers' binary
shifts. We provide here the details of $n$-unitary
shifts generalized by \cite {mC87} M.Choda for the convenience
of the reader.

Let $n$ be a positive integer. We treat a pair of sets $Q$ and
$S$ of integers satisfying the following condition $(*)$ for
some integer $m$:
\[
(*)
\begin {cases}
Q= (i(1), i(2), \cdots, i(m)), &0 \leq i(1) < i(2) < \cdots
< i(m), \\
S= (j(1), j(2), \cdots, j(m)), & j(l) \in {\mathbb N}, \quad j(l) \leq n-1\\
&\text {\ for \ } l= 1, 2, \cdots, m.
\end {cases}
\]
\begin {definition}
A unital $*$-endomorphism $\Psi$ of $R$ is called an $n$-unitary
shift of $R$ if there is a unitary $u \in R$ satisfying the following:\\
(1)$u^n= 1$;\\
(2)$R$ is generated by $\{u, \Psi (u), \Psi^2 (u), \cdots, \}$;\\
(3)$\Psi^k (u) u= u \Psi^k (u)$ or $\Psi^k (u) u= \gamma u \Psi^k
(u)$
for all $k= 1,2,\cdots$, where $\gamma = \exp (2 \pi \sqrt {-1} / n)$.\\
(4) for each $(Q, S)$ satisfying $(*)$, there are an integer $k (\geq
0)$ and a nontrivial $\lambda \in \mathbb {T} = \{ \mu \in \mathbb {C};
| \mu | = 1 \}$ such that
\[
  \Psi^k (u) u(Q, S)= \lambda u(Q, S) \Psi^k (u),
\]
where  $u(Q, S)$ is defined by
\[
  u(Q, S) = {\Psi^{i(1)} (u)}^{j(1)} {\Psi^{i(2)} (u)}^{j(2)}
  \cdots  {\Psi^{i(m)} (u)}^{j(m)}.
\]
\end {definition}

The unitary $u$ is called a generator of $\Psi$. Put $S (\Psi; u)=
\{ k; \Psi^k (u) u= \gamma u \Psi^k (u) \}$. Note that the above
condition (2) gives some rigidity on $S (\Psi; u)$.
The Jones index $[R: \Psi (R)]$ is $n$.

One interesting example of 
\[
  S_1= S (\Psi_1; u_1) = \{ 1, 3, 6, 10, 15, \cdots, 
  \frac {1} {2} i(i+1),  \cdots \}_{i \in {\mathbb N}}, 
\]
which corresponds to the $n$-stream 
$\{ 0 1 0 1 0 0 1 0 0 0 1 0 0 0 0 1 0 0 0 0 0 1 \cdots \}$. 
It is pointed out in \cite {mC87} that the relative commutant
$\Psi_1^k (R)' \cap R$ is always trivial for all $k$! That is, our
question for $A= \mathbb {C}$ is answered by this example.

\section {Main Theorem} \label {S:main}

\begin{theorem}
For any finite dimensional $C^*$-algebra $A$ with any trace vector
$\vec s$ whose components are rational numbers, we give an endomorphism
$\Phi$ of the hyperfinite $II_1$ factor $R$  such that:
\[
  \forall \ k \ \in \ \mathbb {N}, \
  \Phi^k (R)' \cap R= \otimes^k A.
\]
The canonical trace $\tau$ on $R$ extends the trace vector $\vec s$
on $A$.
\end{theorem}

$A$ is characterized by its trace vector $\vec s$ and its dimension
vector $\vec t$.
\begin {align*}
   &\vec s= [\frac{b_1}{c_1}, \ \frac{b_2}{c_2}, \cdots,
   \frac{b_n}{c_n}]\\
   &\vec t= [a_1, \ a_2, \cdots, a_n]\\
   &\frac{b_1}{c_1} a_1+ \frac{b_2}{c_2} a_2+ \cdots + \frac{b_n}{c_n} a_n=
   1\\
   &b_1,\  c_1, \  a_1, \  b_2,\  c_2, \  a_2, \cdots, b_n, \ c_n, \ a_n \in
   {\mathbb N}
\end {align*}

Put $d= c_1 c_2 \cdots c_n$. We can embed $A$ into $B \subset M_{d}
({\mathbb C})$ via
\begin{align*}
  &A \simeq
  M_{a_1} ({\mathbb C}) \oplus M_{a_2} ({\mathbb C})
  \oplus \cdots \oplus M_{a_n} ({\mathbb C}) \subset\\
  &M_{a_1} ({\mathbb C}) \otimes M_{d_1} ({\mathbb C})
   \oplus M_{a_2} ({\mathbb C}) \otimes M_{d_2} ({\mathbb C})
  \oplus \cdots \oplus M_{a_n} ({\mathbb C}) \otimes M_{d_n} ({\mathbb
  C})\\
  &= B \subset M_{d} ({\mathbb C})\\
  &{\rm where \ } {d}_1= \frac{d
  b_1}{c_1}, \ {d}_2= \frac{d b_2}{c_2}, \cdots, {d}_n=
  \frac{d b_n}{c_n}
\end{align*}

For each $i$, $1 \leq i \leq n$, $M_{a_i} (\mathbb {C}) \subset A \subset B$ (the
former being not a unital embedding) is generated by $p_i, q_i \in
\mathcal {U}(\mathbb {C}^{a_i})$ with:
\[
  p_i^{a_i}= q_i^{a_i}= {\bf 1}_{M_{a_i}(\mathbb {C})}; \quad
  \gamma_i= \exp( 2 \pi \sqrt {-1} / a_i), \quad p_i q_i = \gamma_i  q_i p_i
\]
where $p_i$ is the diagonal matrix in $M_{a_i} (\mathbb {C})$,
$[ 1 \ \gamma_i \ \gamma_i^2 \cdots \gamma_i^{a_i-1} ]$,
 and
$q_i$ is the permutation matrix in $M_{a_i} (\mathbb {C})$, $( 1 \ 2
\ 3 \cdots a_i )$.

For each $i$, $1 \leq i \leq n$, $M_{{d}_i} (\mathbb {C}) \subset A' \cap B \subset B$
(the former being not a unital embedding) is generated by ${\mathfrak p}_i,
{\mathfrak q}_i \in \mathcal {U}(\mathbb {C}^{{d}_i})$ with:
\[
  {\mathfrak p}_i^{{d}_i}= {\mathfrak q}_i^{{d}_i}= {\bf 1}_{M_{{d}_i}(\mathbb {C})}; 
  \quad
  \rho_i= \exp( 2 \pi \sqrt {-1} / {d}_i), \quad {\mathfrak p}_i {\mathfrak q}_i = \rho_i  
  {\mathfrak q}_i
  {\mathfrak p}_i
\]
where ${\mathfrak p}_i$ is the diagonal matrix in $M_{{d}_i} (\mathbb {C})$,
$[ 1 \ \rho_i \ \rho_i^2 \cdots \rho_i^{{d}_i-1} ]$,   and
${\mathfrak q}_i$ is the permutation matrix in $M_{{d}_i} (\mathbb {C})$, $( 1 \
2 \ 3 \cdots {d}_i)$.

Define $v \in M_{d} (\mathbb {C})$ to be the permutation matrix:
\[
v= (a_1 {d}_1 \  \ (a_1 {d}_1+ a_2 {d}_2) \ \
\cdots \ \  (a_1 {d}_1+ a_2 {d}_2+ \cdots + a_n {d}_n))
\]
Then $v^n= {\bf 1}$. $B$ and $v$  generate $M_{d} (\mathbb {C})$.

Define $r:= s v s$, while $s$ is a diagonal matrix in $M_{d} ({\mathbb C})$,
\[
  [0 \ 0 \ \cdots \ 0 \ 1_{a_1 {d}_1} \ 0 \ 0 \
  \cdots \ 0 \ 1_{a_1 {d}_1+ a_2 {d}_2} \cdots \ 0 \ 0 \cdots \ 0 \ 1_{a_1 {d}_1+ 
  a_2 {d}_2 + \cdots + a_n {d}_n}]
\]
$s$ lies in $B$ and $r^n= s$. Since $v= {\bf 1}-s+ r$,
$B$ and $r$ generates $M_{d} (\mathbb {C})$.

\begin {lemma}
$<B, r> \simeq M_{d} (\mathbb {C})$ is of the form:
\[
  B+ B r B+ B r^2 B+ \cdots + B r^{n-1} B.
\]
\end {lemma}

\begin {proof}
It suffices to observe that
$sBs$ is abelian and  $Ad v$ sends $sBs$ onto itself.
\begin{align*}
   & [s, v]= [s, r]=0, \quad r= s v s= v s = s v\\
   &r^2= (s v s) (s v s)= s v^2 s, \quad  
   r^* = s v^* s= s v^{n-1} s= r^{n-1}\\
   &r B r= (s v s) B  (s v s) = s v (s B s) v s
   = s v^2 (s B s) s \subset r^2 B\\
   &r^* r= r r^*= r^n= s \in B
\end {align*}
\end {proof}

Define $w= \sum_{i=1}^n \gamma^{i-1} {\bf 1}_{M_{a_i}
(\mathbb {C})} \otimes  {\bf 1}_{M_{{d}_i}
(\mathbb {C})}$, where 
\[
  \gamma= \exp (2 \pi \sqrt {-1} / n),\quad
  w^n = {\bf 1}_{M_{d} (\mathbb {C})}. 
\]
Note that $w$ is in the center of $B$.
Observe that:\\
(1) $Ad w$ acts trivially on $B$.\\
(2) $Ad w (r)= \gamma r$.

Now we  construct a tower of inclusions of finite dimensional $C^*$-
algebras $M_k$ with a trace $\tau$. The ascending union $M= \cup_{k
\in \mathbb {N}} M_k$ contains infinitely many copies of $M_{d} (\mathbb
{C})$ and thus, infinitely many copies of $B$. Number them respectively by 
$r_1, A_1 \subset B_1, w_1$,
$r_2, A_2 \subset B_2, w_2$, $r_3, A_3 \subset B_3, w_3$, $\cdots$.

We endow on the union $M$ the following properties:
\begin{align*}
  &[r_l, B_m]= 0, \ \text{if} \ m < l\\
  &r_l r_m = \gamma r_m r_l \  \text{if} \  l-m \in
  S_1= \{ 1, 3, 6, 10, 15, \cdots \}\\
  &r_l r_m = r_m r_l, \ \text{if \ } l-m \notin S_1
\end{align*}

There is no twist in the relation
between $A_m$ and $A_l$, $m < l$,
where:
\begin{align*}
  &A_m= \oplus_{i=1}^{n} <(p_{i})_m, (q_{i})_m>\\
  &(p_{i})_m (q_{i})_m= \gamma_{i} (q_{i})_m
  (p_{i})_m\\
  &A_l= \oplus_{j=1}^{n} <(p_{j})_l, (q_{j})_l>\\
   &(p_{j})_l (q_{j})_l= \gamma_{j} (q_{j})_l
  (p_{j})_l\\
  &(p_{j})_l (p_{i})_m = (p_i)_m (p_{j})_l\\
   &(p_{j})_l (q_{i})_m= (q_{i})_m (p_{j})_l\\
  &(q_{j})_l (q_{i})_m = (q_{i})_m (q_{j})_l
\end{align*}

Define
\[
  S_2:= \{ \frac{1}{2} i (i+1) \ | \ i= 1 \mod 3 \},
  \quad  S_3:= \{ \frac{1}{2} i (i+1) \ | \ i= 2 \mod 3 \}
\]

We add a twist in the relation between $(A' \cap B)_m$
and $(A' \cap B)_l$, $m < l$ where:
\begin{align*}
  &(A' \cap B)_m= \oplus_{i=1}^{n} <({\mathfrak p}_{i})_m, ({\mathfrak q}_{i})_m>\\
  &({\mathfrak p}_{i})_m ({\mathfrak q}_{i})_m = \rho_i ({\mathfrak q}_{i})_m 
({\mathfrak p}_{i})_m\\
   &(A' \cap B)_l= \oplus_{j=1}^{n} <({\mathfrak p}_{j})_l, ({\mathfrak q}_{j})_l>\\
  &({\mathfrak p}_{j})_l ({\mathfrak q}_{j})_l = \rho_j ({\mathfrak q}_{j})_l 
({\mathfrak p}_{j})_l\\
  &({\mathfrak p}_{j})_l ({\mathfrak p}_{i})_m= ({\mathfrak p}_{i})_m ({\mathfrak p}_{j})_l\\
  &({\mathfrak p}_{j})_l ({\mathfrak q}_{i})_m =({\mathfrak q}_{i})_m ({\mathfrak p}_{j})_l\\
    &{\rm If \ } l-m \in S_2  {\rm \ and \ } l > m^2,\\
  &{\rm then \ } ({\mathfrak q}_{j})_l ({\mathfrak p}_{i})_m=  \rho_i^{- \delta_{i j}}
  ({\mathfrak p}_{i})_m ({\mathfrak q}_{j})_l.\\
   &{\rm If \ } l-m \notin S_2  {\rm \ and \ } l > m^2,\\
  &({\mathfrak q}_{j})_l ({\mathfrak p}_{i})_m=  
  ({\mathfrak p}_{i})_m ({\mathfrak q}_{j})_l. \\
    &{\rm If \ } l-m \in S_3  {\rm \ and \ } l > m^2,\\
   &({\mathfrak q}_{j})_l ({\mathfrak q}_{i})_m=  \rho_i^{\delta_{i j}}
  ({\mathfrak q}_{i})_m ({\mathfrak q}_{j})_l.\\
    &{\rm If \ } l-m \notin S_3  {\rm \ and \ } l > m^2,\\
   &({\mathfrak q}_{j})_l ({\mathfrak q}_{i})_m=  
  ({\mathfrak q}_{i})_m ({\mathfrak q}_{j})_l.
\end{align*}

The construction is an induction process. 
We embed  $M_1 \simeq M_{d} (\mathbb {C})$ into $\otimes^2 M_{d} (\mathbb {C})$,
by sending any element $x \in M_1$ to $x \otimes {\bf 1}$.
$M_1$ is equipped
with the trace $\frac{1}{d} Tr$.

Observe that $2-1 = 1 \in S_2$. 
There is a twist in the relation of $B_1$ and $B_2$, where
$B_2$ is generated by
$({p}_{j})_2$, $({q}_{j})_2$,
$({\mathfrak p}_{j})_2$, and $({\mathfrak q}_{j})_2$, $1 \leq j \leq n$.

Put 
\begin{align*}
   &({p}_{j})_2= {\bf 1} \otimes p_i \in \otimes^2 M_{d} ({\mathbb C})\\
   &({q}_{j})_2= {\bf 1} \otimes q_i \in \otimes^2 M_{d} ({\mathbb C})\\
   &({\mathfrak p}_{j})_2= {\bf 1} \otimes {\mathfrak p}_{j} 
   \in \otimes^2 M_{d} ({\mathbb C})\\
   &({\mathfrak q}_{j})_2= ({\mathfrak q}_{j} + {\bf 1}- {\bf 1}_{M_{d_j} 
({\mathbb C})}) \otimes
   {\mathfrak q}_{j} \in \otimes^2 M_{d} ({\mathbb C})
\end{align*}

Note that ${\mathfrak q}_{j} + {\bf 1}- {\bf 1}_{M_{d_j} ({\mathbb C})}
\in {\mathcal U} ({\mathbb C}^d)$.  We
have:
\begin {align*}
    &({\mathfrak p}_{j})_2^{d_j}= ({\mathfrak q}_{j})_2^{d_j}
    = {\bf 1} \otimes {\bf 1}_{M_{d_j} ({\mathbb C})}\\
   &({\mathfrak p}_{j})_2 ({\mathfrak q}_{j})_2 
   =\rho_j ({\mathfrak q}_{j})_2 ({\mathfrak p}_{j})_2\\
   &({\mathfrak p}_{j})_2 ({\mathfrak p}_{i})_1 =
   ({\mathfrak p}_{i})_1 ({\mathfrak p}_{j})_2\\
   &({\mathfrak p}_{j})_2 ({\mathfrak q}_{i})_2 =
   ({\mathfrak q}_{i})_1 ({\mathfrak p}_{j})_2\\
   &({\mathfrak q}_{j})_2 ({\mathfrak p}_{i})_1 =
    \rho_i^{-\delta_{ij}}
   ({\mathfrak p}_{i})_2 ({\mathfrak q}_{j})_2\\
   &({\mathfrak q}_{j})_2 ({\mathfrak q}_{i})_2 =
   ({\mathfrak q}_{i})_2 ({\mathfrak q}_{j})_2
\end{align*}

Thus
\[
  B_2 \simeq B= \oplus_{i=1}^{n} M_{a_i} ({\mathbb C}) \otimes  
  M_{d_i} ({\mathbb C}) \subset M_{d}  ({\mathbb C})
\]

Observe $2-1=1 \in S_1$. Define $r_2: = w \otimes r$. We have the
following properties:
\begin{align*}
   &[w, {\mathfrak q}_{j} + {\bf 1}- {\bf 1}_{M_{d_j} ({\mathbb C})}]= 0\\
   &<B_2, r_2> \simeq M_{d} (\mathbb {C})\\
   &[B_1, r_2]= 0 \quad r_1 r_2= \gamma r_2 r_1\\
   &M_2= <B_1, r_1, B_2, r_2> = \otimes^2 M_{d} ({\mathbb C})
\end{align*}
$M_2$ is equipped with a unique normalized trace $\tau$.

Assume that we have obtained $M_k= <B_1, r_1, B_2, r_2, \cdots, B_k,
r_k>$ isomorphic to $ \otimes^k M_{d} ({\mathbb C})$ with the trace $\tau$.
We embed $M_k$ into $M_k \otimes {M_{d} (\mathbb {C})}$ by sending
$x \in M_k$ to $x \otimes {\bf 1}_{M_d (\mathbb {C})}$.

Define $B_{k+1}$ by its generators: $({p}_{j})_{k+1}$,
$({q}_{j})_{k+1}$, $({\mathfrak p}_{j})_{k+1}$,
$({\mathfrak q}_{j})_{k+1}$, where $1 \leq j \leq n$:
\begin{align*}
   &({p}_{j})_{k+1} := {\bf 1} \otimes {\bf 1} \otimes \cdots \otimes {\bf 1} \otimes p_j\\
   &({q}_{j})_{k+1} := {\bf 1} \otimes {\bf 1} \otimes \cdots \otimes {\bf 1} \otimes q_j\\
   &({\mathfrak p}_{j})_{k+1}:= {\bf 1} \otimes {\bf 1} \otimes \cdots \otimes {\bf 1} \otimes 
    {\mathfrak p}_{j}\\
   &({\mathfrak q}_{j})_{k+1}:= \\
   &[({\mathfrak q}_{j} + {\bf 1}- {\bf 1}_{M_{d_j} ({\mathbb C})})^{\beta_1})
    \otimes \cdots
    \otimes ({\mathfrak q}_{j} + {\bf 1}- {\bf 1}_{M_{d_j} ({\mathbb C})})^{\beta_k}
    ] \cdot \\
   &[({\mathfrak p}_{j} + {\bf 1}- {\bf 1}_{M_{d_j} ({\mathbb C})})^{\beta_{k+1}})
    \otimes \cdots
    \otimes ({\mathfrak p}_{j} + {\bf 1}- {\bf 1}_{M_{d_j} ({\mathbb C})})^{\beta_{2k}}
    ] \otimes {\mathfrak q}_{j}\\
    &{\rm Where \ } 1 \leq  i \leq k:\\
    &\beta_i= 1, \quad {\rm if \ } k+1-i \in S_2; \quad
   \beta_l= 0, \quad {\rm if \ } k+1-i \notin S_2\\
      &\beta_{k+i}= 1, \quad {\rm if \ } k+1-i \in S_3; \quad
   \beta_{k+i}= 0, \quad {\rm if \ } k+i-1 \notin S_3
\end{align*}

We have:
\begin{align*}
   &({p}_{j})_{k+1}^{a_j}= ({q}_{j})_{k+1}^{a_j}= 
   \otimes^k {\bf 1} \otimes {\bf 1}_{M_{a_j}}({\mathbb C})\\
   &({p}_{j})_{k+1} ({q}_{j})_{k+1}= \gamma_j
   ({q}_{j})_{k+1} ({p}_{j})_{k+1}\\
   &({\mathfrak p}_{j})_{k+1}^{d_j}= ({\mathfrak q}_{j})_{k+1}^{d_j}= 
   \otimes^k {\bf 1} \otimes {\bf 1}_{M_{d_j}}({\mathbb C})\\
   &({\mathfrak p}_{j})_{k+1} ({\mathfrak q}_{j})_{k+1}= \rho_j
   ({\mathfrak q}_{j})_{k+1} ({\mathfrak p}_{j})_{k+1}
\end{align*}
Therefore $B_{k+1}$ is isomorphic to $B$.

The commutation relations are given below.
\begin{align*}
   &l < k+1:\\ 
   &[({p}_{j})_{k+1}, ({p}_{i})_{l}]= 0 \\
   &[({p}_{j})_{k+1}, ({q}_{i})_{l}]= 0 \\
   &[({q}_{j})_{k+1}, ({p}_{i})_{l}]= 0 \\
   &[({q}_{j})_{k+1}, ({q}_{i})_{l}]= 0 \\
   &A_{k+1} \cdot A_l= A_l \cdot A_{k+1}\\
    &[({\mathfrak p}_{j})_{k+1}, ({\mathfrak p}_{i})_{l}]= 0 \\
     &[({\mathfrak p}_{j})_{k+1}, ({\mathfrak q}_{i})_{l}]= 0 
\end{align*}

The anti-commutation relations are given below.
\begin{align*}
   &{\rm If \ } k+1-l \in S_2 {\rm \ and \ }k+1 > l^2,\\ 
     &({\mathfrak q}_{j})_{k+1} ({\mathfrak p}_{i})_{l}= \rho_i^{-\delta_{ij}}
    ({\mathfrak p}_{i})_{l} ({\mathfrak q}_{j})_{k+1}.\\
    &{\rm If \ } k+1-l \notin S_2 {\rm \ and \ }k+1 > l^2,\\
     &({\mathfrak q}_{j})_{k+1} ({\mathfrak p}_{i})_{l}= 
    ({\mathfrak p}_{i})_{l} ({\mathfrak q}_{j})_{k+1}.\\
   &{\rm If \ } k+1-l \in S_3 {\rm \ and \ }k+1 > l^2, \\
     &({\mathfrak q}_{j})_{k+1} ({\mathfrak q}_{i})_{l}= \rho_i^{\delta_{ij}}
    ({\mathfrak q}_{i})_{l} ({\mathfrak q}_{j})_{k+1}.\\
    &{\rm If \ } k+1-l \notin S_3 {\rm \ and \ }k+1 > l^2,\\
    &({\mathfrak q}_{j})_{k+1} ({\mathfrak q}_{i})_{l}= 
    ({\mathfrak q}_{i})_{l} ({\mathfrak q}_{j})_{k+1}.
\end{align*}

Define
\begin {align*}
   &r_{k+1}:= w^{\alpha_1} \otimes w^{\alpha_2} \otimes \cdots \otimes w^{\alpha_k}
   \otimes r\\
   &\alpha_i= 1, \quad {\rm if \ } k+1-i \in S_1; \quad
   \alpha_i= 0, \quad {\rm otherwise.}
\end{align*}

We have the following properties:
\begin{align*}
   &[w, {\mathfrak q}_{j} + {\bf 1}- {\bf 1}_{M_{d_j} ({\mathbb C})}]= [w, 
   {\mathfrak p}_{j} + {\bf 1}- {\bf 1}_{M_{d_j} ({\mathbb C})}]= 0\\
   &<B_{k+1}, r_{k+1}> \simeq M_n (\mathbb {C})\\
   &[B_j, r_{k+1}]= 0 \quad 1 \leq j \leq k\\
    &r_{k+1}  r_j= \gamma r_j r_{k+1} \quad {\rm if \ } k+1-j \in S_1\\
    &r_{k+1}  r_j= r_j r_{k+1} \quad {\rm if \ } k+1-j \notin S_1\\
   &M_{k+1}= <M_k, B_{k+1}, r_{k+1}> = \otimes^{k+1} M_{d} ({\mathbb C})
\end{align*}
There is  a unique normalized trace $\tau$ on $M_{k+1}$.

By induction we have constructed the ascending tower of inclusions of
finite dimensional $C^*$-algebras with the desired properties.

We now explore some useful properties of the finite dimensional
$C^*$-algebra, $M_k$.

\begin {lemma}
For all $k$, $M_k$ is the linear span of the words, $x_1 \cdot x_2
\cdot x_3  \cdots  x_k$, where $x_j \in (M_{d} (\mathbb {C}))_j= <B_j,
r_j>$.
\end {lemma}

\begin {proof}
It suffices to prove $x_j \cdot x_i$ is in $M_i \cdot <B_j, r_j>= M_i \cdot 
(M_{d} (\mathbb {C}))_j$, where $i < j$ and $M_i= \otimes^i M_d({\mathbb C})$.

\begin{align*}
&r_j \cdot B_i= B_i \cdot r_j\\
&{\rm Either \quad } r_j \cdot r_i= r_i \cdot r_j \quad
{\rm  or \quad } r_j \cdot r_i= \gamma r_i \cdot r_j\\
&A_j \cdot A_i= A_i \cdot A_j \subset M_i \cdot A_j\\
&A_j \cdot (A' \cap B)_i= (A' \cap B)_i \cdot A_j \subset M_i \cdot A_j\\
&(A' \cap B)_j \cdot A_i= A_i \cdot (A' \cap B)_j \subset M_i \cdot A_j\\
&{\rm Note \ that:}\\ 
&Ad ({\mathfrak p}_{\iota} + {\bf 1}- {\bf 1}_{M_{d_\iota} ({\mathbb C})}) 
({\mathfrak q}_{\iota})= \rho_{\iota} {\mathfrak q}_{\iota}\\
&Ad ({\mathfrak q}_{\iota} + {\bf 1}- {\bf 1}_{M_{d_\iota} ({\mathbb C})}) 
({\mathfrak p}_{\iota})= \rho_{\iota}^{-1} {\mathfrak p}_{\iota}\\
&Ad ({\mathfrak p}_{\iota} + {\bf 1}- {\bf 1}_{M_{d_\iota} ({\mathbb C})}) 
({\mathfrak q}_{\iota}+  {\bf 1}- {\bf 1}_{M_{d_\iota} ({\mathbb C})})= 
\rho_{\iota} {\mathfrak q}_{\iota}+  {\bf 1}- {\bf 1}_{M_{d_\iota} ({\mathbb C})}\\
&Ad ({\mathfrak q}_{\iota} + {\bf 1}- {\bf 1}_{M_{d_\iota} ({\mathbb C})}) 
({\mathfrak p}_{\iota} + {\bf 1}- {\bf 1}_{M_{d_\iota} ({\mathbb C})})= 
 \rho_{\iota}^{-1} {\mathfrak p}_{\iota} + {\bf 1}- {\bf 1}_{M_{d_\iota} ({\mathbb C})}\\
&(A' \cap B)_j \cdot (A' \cap B)_i \subset M_i \cdot (A' \cap B)_j\\
&{\rm In \ short, \ } B_j \cdot B_i \subset M_i \cdot B_j.
\end{align*}

\begin{align*}
&{\rm Note \ that:}\\ 
&Ad ({\mathfrak p}_{\iota} + {\bf 1}- {\bf 1}_{M_{d_\iota} ({\mathbb C})}) 
(w)= w\\
&Ad ({\mathfrak q}_{\iota} + {\bf 1}- {\bf 1}_{M_{d_\iota} ({\mathbb C})}) 
(w)= w\\
&Ad ({\mathfrak p}_{\iota} + {\bf 1}- {\bf 1}_{M_{d_\iota} ({\mathbb C})}) 
(r) \in M_d ({\mathbb C})\\
&Ad ({\mathfrak q}_{\iota} + {\bf 1}- {\bf 1}_{M_{d_\iota} ({\mathbb C})}) 
(r) \in M_d ({\mathbb C})\\
&{\rm In \ short, \ } B_j \cdot r_i \subset M_i \cdot B_j.\\
\end{align*}
\end {proof}

\begin {lemma}
Consider the GNS-construction of the pair $(M, \tau)$ described above.
The weak closure of 
$M$ is the hyperfinite $II_1$ factor.
\end {lemma}

\begin{proof}
There is a unique tracial state on $M_k= \otimes^k M_d ({\mathbb C})$ for all $k \in
{\mathbb N}$, and hence a unique tracial state on $M$, a $d^{\infty}$ UHF-algebra. 
\end {proof}

Define a unital *-endomorphism, $\Phi$, on $R$ by
sending $B_k$ to $B_{k+1}$ and sending $r_k$ to $r_{k+1}$: 
\begin{align*}
   &1 \leq i \leq n:\\
   &\Phi ((p_i)_k)= (p_i)_{k+1}\\
    &\Phi ((q_i)_k)= (p_i)_{k+1}\\
   &\Phi (({\mathfrak p}_i)_k)= ({\mathfrak p}_i)_{k+1}\\
    &\Phi (({\mathfrak q}_i)_k)= ({\mathfrak q}_i)_{k+1}\\
    &\Phi (r_k)= r_{k+1}\\
\end{align*}

We observe that $\Phi (R)$ is a hyperfinite $II_1$ factor and
\[
  [R: \Phi (R)]= d^2 < \infty.
\]

\begin {lemma}
The relative commutant $\Phi^k (R)' \cap R$ is exactly
$\otimes_{i=1}^k A$, on which the trace of $R$ is the product trace
given by the vector $\vec s$.
\end {lemma}

\begin {proof}
Because of our decomposition in Lemma 2 and Lemma 3, $R$ can be
written as
\[
  (\sum_{i=0}^{n-1} B_1 r_1^{i} B_1) \cdot (\sum_{i=0}^{n-1}
  B_2 r_2^{i} B_2) \cdots  (\sum_{i=0}^{n-1} B_k r_k^{i} B_k)
  \cdot \Phi^k (R).
\]

Assume $x \in \Phi^k (R)' \cap R$.
Let $\vec {\alpha}= (i_1, i_2, \cdots, i_k)$ be a multi-index. 
\begin{align*}
   &0 \leq i_1, i_2, \cdots, i_k \leq n-1\\
   & x= \sum_{\vec {\alpha} \in \{ 0, \ 1,  \ \cdots, n-1 \}^k} y_1^{\vec {\alpha}} 
  r_1^{i_1} z_1^{\vec{\alpha}}  y_2^{\vec {\alpha}} r_2^{i_2} z_2^{\vec {\alpha}} 
  \cdots y_k^{\vec {\alpha}} r_k^{i_k} z_k^{\vec {\alpha}}
  \cdot y^{\vec {\alpha}}\\
   &y_1^{\vec {\alpha}}, z_1^{\vec {\alpha}} \in B_1 \\
   &y_2^{\vec {\alpha}}, z_2^{\vec {\alpha}} \in B_2\\
   &\cdots\\ 
   &y_k^{\vec{\alpha}}, z_k^{\vec {\alpha}} \in B_k\\
   &y^{\vec {\alpha}} \in \Phi^k (R)
\end{align*}
Note that $\Phi^k (R)$ is the weak closure of $\{ \Phi^k (M_j) \}_{j=1}^{\infty}$.

For every $\epsilon > 0$,  there exists an integer $j \in {\mathbb
N}$ such that
\begin{align*}
   &\forall {\vec {\alpha}} \quad \exists z^{\vec {\alpha}} \in \Phi^k (M_j)
   \subset <B_{k+1}, r_{k+1}, \cdots, B_{k+j}, r_{k+j}>\\
   &\| x-\sum_{\vec {\alpha} \in \{ 0, \ 1,  \ \cdots, n-1 \}^k} 
   y_1^{\vec {\alpha}} r_1^{j_1} z_1^{\vec
  {\alpha}}
  y_2^{\vec {\alpha}} r_2^{j_2} z_2^{\vec {\alpha}} \cdots y_k^{\vec {\alpha}} r_k^{j_k}
  z_k^{\vec {\alpha}}
  \cdot z^{\vec {\alpha}} \|_{2, \tau} < \delta\\
  &\delta= (\sqrt{\frac{n}{d}})^k \epsilon
\end{align*}

Put $L= l( l+1)/2+1$ for some integer $l > 2 (k+j)$ and $l= 0 \mod 3$.
We have the following properties:
\begin {align*}
&[r_L, B_1]= [r_L, B_2]= \cdots = [r_L, B_{k+j}]= 0\\
&[r_L, r_2]= [r_L, r_3]= \cdots = [r_L, r_{k+j}]= 0\\
&r_L r_1 = \gamma r_1 r_L, \quad r_L r_1^{j_1}
   = \gamma^{j_1} r_1 r_L\\
&r_L r_L^*= r_L^* r_L= s_L\\
&r_L r_1 r_L^*= \gamma r_1 s_L, \quad r_L r_1^{j_1} r_L^*= \gamma^{j_1} r_1 s_L\\
&{\rm for \ } 0 \leq m \leq {n-1}, \quad r_L^m r_1^{g_1} (r_L^*)^m= \gamma^{j_1 m} r_1 s_l\\
&[s_L, B_1]= [s_L, B_2]= \cdots = [s_L, B_{k+j}]= 0\\
&[s_L, r_1]= [s_L, r_2]= \cdots = [s_L, r_{k+j}]= 0
\end {align*}

Therefore we claim:
\begin {align*}
&\| (x-\sum_{\vec {\alpha} \in \{ 0, \ 1,  \ \cdots, n-1 \}^k}
y_1^{\vec {\alpha}} r_1^{j_1} z_1^{\vec
  {\alpha}}
  y_2^{\vec {\alpha}} r_2^{j_2} z_2^{\vec {\alpha}} \cdots y_k^{\vec {\alpha}} r_k^{j_k}
  z_k^{\vec {\alpha}}
  \cdot z^{\vec {\alpha}}) s_L \|_{2, \tau}=\\
&\| (x-\sum_{\vec {\alpha} \in \{ 0, \ 1,  \ \cdots, n-1 \}^k}
y_1^{\vec {\alpha}} r_1^{j_1} z_1^{\vec
  {\alpha}}
  y_2^{\vec {\alpha}} r_2^{j_2} z_2^{\vec {\alpha}} \cdots y_k^{\vec {\alpha}} r_k^{j_k}
  z_k^{\vec {\alpha}}
  \cdot z^{\vec {\alpha}}) \frac{1}{n} \sum_{m=0}^{n-1} r_L^m {r_L^*}^m \|_{2, \tau}=\\
&\frac{1} {n} \| \sum_{\vec {\alpha}} \sum_m (r_L^m x {r_L^*}^m-
y_1^{\vec {\alpha}} r_L^m r_1^{j_1} (r_L^*)^m z_1^{\vec
  {\alpha}}
  y_2^{\vec {\alpha}} r_2^{j_2} z_2^{\vec {\alpha}} \cdots y_k^{\vec {\alpha}} j_1^{j_k}
  z_k^{\vec {\alpha}}
  \cdot z^{\vec {\alpha}}) \|_{2, \tau}=\\
&\frac{1} {n} \| \sum_{\vec {\alpha}} \sum_m (x - y_1^{\vec
{\alpha}}  r_1^{j_1 m} z_1^{\vec
  {\alpha}}
  y_2^{\vec {\alpha}} r_2^{j_2} z_2^{\vec {\alpha}} \cdots y_k^{\vec {\alpha}} r_k^{j_k}
  z_k^{\vec {\alpha}}
  \cdot z^{\vec {\alpha}}) s_L \|_{2, \tau}=\\
&\| (x-\sum_{\vec {\alpha} \in \{ 0, \ 1, \cdots, j-1 \}^k, j_1= 0}
y_1^{\vec {\alpha}} z_1^{\vec
  {\alpha}}
  y_2^{\vec {\alpha}} r_2^{j_2} z_2^{\vec {\alpha}} \cdots y_k^{\vec {\alpha}} r_k^{j_k}
  z_k^{\vec {\alpha}}
  \cdot z^{\vec {\alpha}}) s_L \|_{2, \tau}=\\
&\sqrt{\frac{j}{n}} \| x-\sum_{\vec {\alpha} \in \{ 0, \ 1, \cdots,
n-1 \}^k, g_1= 0} y_1^{\vec {\alpha}} z_1^{\vec
  {\alpha}}
  y_2^{\vec {\alpha}} r_2^{j_2} z_2^{\vec {\alpha}} \cdots y_k^{\vec {\alpha}} r_k^{j_k}
  z_k^{\vec {\alpha}}
  \cdot z^{\vec {\alpha}}  \|_{2, \tau}
\end{align*}

Since
\begin{align*}
&\{ x, \ y_1^{\vec {\alpha}}, \ z_1^{\vec
  {\alpha}}, \ y_2^{\vec {\alpha}}, \ r_2^{j_2}, \ z_2^{\vec {\alpha}}, \cdots,
 y_k^{\vec{\alpha}},
  \ r_k^{j_k}, \ z_k^{\vec {\alpha}}, \ z^{\vec {\alpha}} \} \subset \{ s_L, \ r_L, \ B_L
  \}'\\
&{\rm and \ } \tau (s_L)= \frac{n}{d}.
\end {align*}
Note that $\{ s_L, \ r_L, \ B_L \}'' = M_d({\mathbb C})_L$ is a type
$I$ factor \cite{sP83}.

By induction,
\begin{align*}
&\| x-\sum_{\vec {\alpha} \in \{ 0, \ 1, \cdots, n-1 \}^k, \ j_1= 0}
y_1^{\vec {\alpha}} z_1^{\vec
  {\alpha}}
  y_2^{\vec {\alpha}} r_2^{j_2} z_2^{\vec {\alpha}} \cdots y_k^{\vec {\alpha}} r_k^{j_k}
  z_k^{\vec {\alpha}}
  \cdot z^{\vec {\alpha}}  \|_{2, \tau} < \sqrt{\frac{d}{n}} \delta\\
&\| x-\sum_{\vec {\alpha} \in \{ 0, \ 1, \cdots, n-1 \}^k, \ j_1=
j_2= 0} y_1^{\vec {\alpha}} z_1^{\vec
  {\alpha}}
  y_2^{\vec {\alpha}} z_2^{\vec {\alpha}} \cdots y_k^{\vec {\alpha}} r_k^{j_k}
  z_k^{\vec {\alpha}}
  \cdot z^{\vec {\alpha}}  \|_{2, \tau} < (\sqrt{\frac{d}{n}})^2 \delta\\
&\cdots \\
&\| x-\sum_{\vec {\alpha} \in \{ 0 \}^k} y_1^{\vec {\alpha}}
z_1^{\vec
  {\alpha}}
  y_2^{\vec {\alpha}} z_2^{\vec {\alpha}} \cdots y_k^{\vec {\alpha}}
  z_k^{\vec {\alpha}}
  \cdot z^{\vec {\alpha}}  \|_{2, \tau} < (\sqrt{\frac{d}{n}})^k \delta= \epsilon\\
\end{align*}

Put $L_1= l_1( l_1+1)/2+1$ for some integer $l_1 > 2 (k+ j)$ and $l_1= 1
\mod 3$. We have the following properties:
\begin{align*}
   &U_{L_1}:= \sum_{m_1=0}^{n -1} ({\mathfrak q}_{m_1})_ {L_1}\\
   &U_{L_1}^d= {\bf 1}\\
   &[U_{L_1}, \Phi (M_{k+j-1})]= 0
\end{align*}

Similarly, put $L_2= l_2( l_2+1)/2+1$ for some integer $l_2 > 2 (k+ j)$
and $l_2= 2 \mod 3$. We have the following properties:
\begin{align*}
   &U_{L_2}:= \sum_{m_2=0}^{n-1} ({\mathfrak q}_{m_2})_ {L_2}\\
   &U_{L_2}^d= {\bf 1}\\
   &[U_{L_2}, \Phi (M_{k+j-1})]= 0
\end{align*}

Therefore
\begin{align*}
&\| x-\sum_{\vec {\alpha} \in \{ 0 \}^k} y_1^{\vec {\alpha}}
z_1^{\vec
  {\alpha}}
  y_2^{\vec {\alpha}} z_2^{\vec {\alpha}} \cdots y_k^{\vec {\alpha}}
  z_k^{\vec {\alpha}}
  \cdot z^{\vec {\alpha}}  \|_{2, \tau} =\\
&\| x-\sum_{\vec {\alpha} \in \{ 0 \}^k} y_1^{\vec {\alpha}}
z_1^{\vec
  {\alpha}}
  y_2^{\vec {\alpha}} z_2^{\vec {\alpha}} \cdots y_k^{\vec {\alpha}}
  z_k^{\vec {\alpha}}
  \cdot z^{\vec {\alpha}}  \frac {1}{d} \sum_{m_3=0}^{d-1} U_{L_1}^{m_3} 
  {U_{L_1}^*}^{m_3} \|_{2,
  \tau}=\\
&\| x-\sum_{\vec {\alpha} \in \{ 0 \}^k} \frac {1}{d}
\sum_{m_3=0}^{d-1} U_{L_1}^{m_3} (y_1^{\vec {\alpha}} z_1^{\vec
  {\alpha}}) {U_{L_1}^*}^{m_3}
\cdots y_k^{\vec {\alpha}}
  z_k^{\vec {\alpha}}
  \cdot z^{\vec {\alpha}}  \|_{2, \tau}=\\
&\| x-\sum_{\vec {\alpha} \in \{ 0 \}^k} \frac {1}{d^2}
\sum_{m_3, m_4=0}^{d-1} U_{L_1}^{m_3} (y_1^{\vec {\alpha}} z_1^{\vec
  {\alpha}}) {U_{L_1}^*}^{m_3}
\cdots y_k^{\vec {\alpha}}
  z_k^{\vec {\alpha}}
  \cdot z^{\vec {\alpha}}  U_{L_2}^{m_4} {U_{L_2}^*}^{m_4}
  \|_{2, \tau}=\\
&\| x-\sum_{\vec {\alpha} \in \{ 0 \}^k} \frac {1}{d^2} \sum_{m_3,
m_4=0}^{d-1} U_{L_2}^{m_4} U_{L_1}^{m_3} (y_1^{\vec {\alpha}}
z_1^{\vec
  {\alpha}}) {U_{L_1}^*}^{m_3} {U_{L_2}^*}^{m_4}
\cdots y_k^{\vec {\alpha}}
  z_k^{\vec {\alpha}}
  \cdot z^{\vec {\alpha}}
  \|_{2, \tau}= \\
&\| x-\sum_{\vec {\alpha} \in \{ 0 \}^k} x_1^{\vec {\alpha}}
  y_2^{\vec {\alpha}} z_2^{\vec {\alpha}} \cdots y_k^{\vec {\alpha}}
  z_k^{\vec {\alpha}}
  \cdot z^{\vec {\alpha}}
  \|_{2, \tau}
\end{align*}

Observe that
\[
  x_1^{\vec {\alpha}}:= \frac {1}{d^2} \sum_{m_3, m_4=0}^{d-1} 
U_{L_2}^{m_4} U_{L_1}^{m_3} 
(y_1^{\vec {\alpha}} z_1^{\vec
  {\alpha}}) {U_{L_1}^*}^{m_3} {U_{L_2}^*}^{m_4}
\]
is the trace-preserving conditional expectation of $y_1^{\vec {\alpha}} z_1^{\vec
  {\alpha}} \in B_1$ onto $A_1$.

By induction,
\begin{align*}
&\| x-\sum_{\vec {\alpha} \in \{ 0 \}^k} x_1^{\vec {\alpha}}
  (y_2^{\vec {\alpha}} z_2^{\vec {\alpha}}) \cdots y_k^{\vec {\alpha}}
  z_k^{\vec {\alpha}}
  \cdot z^{\vec {\alpha}}
  \|_{2, \tau}=\\
&\| x-\sum_{\vec {\alpha} \in \{ 0 \}^k} x_1^{\vec {\alpha}}
  x_2^{\vec {\alpha}}  (y_3^{\vec {\alpha}}
  z_3^{\vec {\alpha}}) \cdots y_k^{\vec {\alpha}}
  z_k^{\vec {\alpha}}
  \cdot z^{\vec {\alpha}}
  \|_{2, \tau}=\\
&\cdots\\
&\| x-\sum_{\vec {\alpha} \in \{ 0 \}^k} x_1^{\vec {\alpha}}
  x_2^{\vec {\alpha}}  \cdots x_k^{\vec {\alpha}}
  \cdot z^{\vec {\alpha}}
  \|_{2, \tau} < \epsilon
\end{align*}
where $x_1^{\vec {\alpha}} \in A_1$, $x_2^{\vec {\alpha}} \in A_2$,
$\cdots$, $x_k^{\vec {\alpha}} \in A_k$.

Note that the von Neumann algebra  $\{ x, \ x_1^{\vec {\alpha}}, \
x_2^{\vec {\alpha}}, \ x_3^{\vec {\alpha}}, \cdots x_k^{\vec
{\alpha}} \}''$ commutes with $\Phi^k (M)$, which is a $II_1$
factor. Any element in the former von Neumann algebra has a scalar
conditional expectation onto $\Phi^k (M)$. In short, the former von
Neumann algebra and $\Phi^k (M)$ are mutually orthogonal \cite {sP83}.

According to the Cauchy-Schwartz inequality
\[
  \| z^{\vec {\alpha}} \|_{2, \tau}^2 \geq
  | \tau (z^{\vec {\alpha}}) |^2
\]
we have:
\begin{align*}
  &\| x-\sum_{\vec {\alpha} \in \{ 0 \}^k} x_1^{\vec {\alpha}}
  x_2^{\vec {\alpha}} x_3^{\vec {\alpha}} \cdots x_k^{\vec {\alpha}}
  \cdot \tau (z^{\vec {\alpha}})  \|_{2, \tau} < \epsilon\\
  &\sum_{\vec {\alpha} \in \{ 0 \}^k} x_1^{\vec {\alpha}}
  x_2^{\vec {\alpha}} x_3^{\vec {\alpha}} \cdots x_k^{\vec {\alpha}}
  \cdot \tau (z^{\vec {\alpha}})  \in A_1 \cdot A_2 \cdots A_k= \otimes^k
  A
\end{align*}
\end{proof}

\section {Application} \label {S:appli}
The Temperley-Lieb algebra \cite {vJ} is generated by projections
$e_i$, $i \in {\mathbb N}$ such that:
\begin{align*}
   &[e_i, e_j]= 0 \quad {\rm if \ } |i -j| \geq 2\\
   &e_i e_{i \pm 1} e_i= \lambda e_i
\end{align*}

We are interested in the case that $\lambda$ is a rational number
and $\lambda^{-1} > 4$.
Put
\[
  \lambda= \frac{p}{q} \quad p, q \in {\mathbb N}
\]

Take $m \in {\mathbb N}$ and $m \geq 3$.
Define  
\[
  A_{1, m}:= vN \{ e_3, e_4, \cdots, e_m \} 
  \quad  A_{0, m}:= vN \{ e_2, e_3, e_4, \cdots, e_m \} 
\]
It is known \cite{vJ} that
\[
  A_{1, m} \subset A_{0, m} \subset M_{q^{[\frac{m}{2}]}} ({\mathbb C})
\]

By the main theorem,
we have an endomorphism $\Phi$ of the hyperfinite $II_1$ factor $R$ such that
\[
  \Phi (R)' \cap R= A_{1, m}
\]

Note that $Q:= \Phi (R)$ is the hyperfinite $II_1$ factor.
Define
\[
   {\mathfrak M}_1= Q \quad {\rm and \quad}
   {\mathfrak M}_l:= vN \{ Q, e_2, e_3, \cdots, e_l \} \quad 2 \leq l \leq m
\]

A corollary to the main theorem is:
\[
  {\mathfrak M}_l' \cap {\mathfrak M}_m = \{ e_2, e_3, \cdots, e_l \}' 
  \cap A_{1, m}= vN \{ e_{l+2}, e_{l+3}, \cdots, e_m \}
\] 
Though ${\mathfrak M}_{m-1} \subset {\mathfrak M}_{m}$ is an irreducible
inclusion of hyperfinite $II_1$ factors,
we do not know a priori whether ${\mathfrak M}_{m-1} \not= {\mathfrak M}_{m}$.
Nor do we know in general whether ${\mathfrak M}_l$ is a factor.

\begin{lemma}
For an element $x \in Q$, 
\[
  \tau (e_2 x)= \lambda \tau (x)
\]
\end{lemma}
\begin{proof}
For every $\epsilon > 0$, we can find $x_k$ in $M_k= \otimes^k M_{q^{[\frac{m}{2}]}} 
({\mathbb C})$ such that $\| x- \Phi (x_k) \|_{2, \tau} < \epsilon$.
Let $\vec {\alpha}= (i_1, i_2, \cdots, i_k)$ be a multi-index. 
\begin{align*}
   &0 \leq i_1, i_2, \cdots, i_k \leq n-1\\
   &x_k= \sum_{\vec {\alpha} \in \{ 0, \ 1,  \ \cdots, n-1 \}^k} y_1^{\vec {\alpha}} 
  r_1^{i_1} z_1^{\vec{\alpha}}  y_2^{\vec {\alpha}} r_2^{i_2} z_2^{\vec {\alpha}} 
  \cdots y_k^{\vec {\alpha}} r_k^{i_k} z_k^{\vec {\alpha}}\\
   &y_1^{\vec {\alpha}}, z_1^{\vec {\alpha}} \in B_1 \quad
   y_2^{\vec {\alpha}}, z_2^{\vec {\alpha}} \in B_2 \quad
   \cdots \quad 
   y_k^{\vec{\alpha}}, z_k^{\vec {\alpha}} \in B_k
\end{align*}

The right-side equation is:
\begin{align*}
  &\tau (\Phi (x_k))= \tau (x_k) \\
  &=\tau (\sum_{\vec {\alpha} \in \{ 0, \ 1,  \ \cdots, n-1 \}^k} y_1^{\vec {\alpha}} 
  r_1^{i_1} z_1^{\vec{\alpha}}  y_2^{\vec {\alpha}} r_2^{i_2} z_2^{\vec {\alpha}} 
  \cdots y_k^{\vec {\alpha}} r_k^{i_k} z_k^{\vec {\alpha}})\\
  &=\tau (\sum_{\vec {\alpha} \in \{ 0 \}^k} y_1^{\vec {\alpha}} 
  z_1^{\vec{\alpha}}  y_2^{\vec {\alpha}} z_2^{\vec {\alpha}} 
  \cdots y_k^{\vec {\alpha}} z_k^{\vec {\alpha}})\\
 &= \tau (\sum_{\vec {\alpha} \in \{ 0 \}^k} \frac {1}{d^2} \sum_{m_3, m_4=0}^{d-1} 
U_{L_2}^{m_4} 
  U_{L_1}^{m_3}   (y_1^{\vec {\alpha}} z_1^{\vec{\alpha}}) {U_{L_1}^*}^{m_3} 
   {U_{L_2}^*}^{m_4}  y_2^{\vec {\alpha}} z_2^{\vec {\alpha}} 
  \cdots y_k^{\vec {\alpha}} z_k^{\vec {\alpha}})\\
   &= \tau (\sum_{\vec {\alpha} \in \{ 0 \}^k} x_1^{\vec {\alpha}} y_2^{\vec {\alpha}} 
   z_2^{\vec {\alpha}} \cdots y_k^{\vec {\alpha}} z_k^{\vec {\alpha}})\\
\end{align*}

Where 
\begin{align*}
   &L_1= l_1( l_1+1)/2+1 \quad l_1 > 2 (k+ 1) \quad 
   l_1= 1 \mod 3\\ 
   &U_{L_1}:= \sum_{m_1=0}^{n -1} ({\mathfrak q}_{m_1})_ {L_1}\\
   &[U_{L_1}, \Phi^2 (M_{k-1})]= 0\\
   &[\Phi(U_{L_1}), e_2]= 0\\
    &L_2= l_2( l_2+1)/2+1 \quad l_2 > 2 (k+ 1) \quad 
    l_2= 2 \mod 3\\
    &U_{L_2}:= \sum_{m_2=0}^{n-1} ({\mathfrak q}_{m_2})_ {L_2}\\
     &[U_{L_2}, \Phi^2 (M_{k-1})]= 0\\
    &[\Phi(U_{L_2}), e_2]=0    
\end{align*}

Observe that
\[
  x_1^{\vec {\alpha}}:= \frac {1}{d^2} \sum_{m_3, m_4=0}^{d-1} 
U_{L_2}^{m_4} U_{L_1}^{m_3} 
(y_1^{\vec {\alpha}} z_1^{\vec
  {\alpha}}) {U_{L_1}^*}^{m_3} {U_{L_2}^*}^{m_4}
\]
is the trace-preserving conditional expectation of $y_1^{\vec {\alpha}} z_1^{\vec
  {\alpha}} \in B_1$ onto $A_1$.

By induction,
\[
  \tau (\Phi (x_k))= \tau (\sum_{\vec {\alpha} \in \{ 0 \}^k} x_1^{\vec {\alpha}}
  x_2^{\vec {\alpha}} \cdots x_k^{\vec {\alpha}})= 
  \tau (\sum_{\vec {\alpha} \in \{ 0 \}^k} \Phi(x_1^{\vec {\alpha}})
  \cdot \Phi(x_2^{\vec {\alpha}}) \cdots \Phi(x_k^{\vec {\alpha}}))
\]
where 
$x_j^{\vec {\alpha}}$ is the trace-preserving conditional expectation of 
$y_j^{\vec {\alpha}} z_j^{\vec
  {\alpha}} \in B_j$ onto $A_j$.

The left-side equation is:
\begin{align*}
  &\tau (e_2 \Phi (x_k))\\
  &=\tau (e_2 \Phi (\sum_{\vec {\alpha} \in \{ 0, \ 1,  \ \cdots, n-1 \}^k} 
  y_1^{\vec {\alpha}} 
  r_1^{i_1} z_1^{\vec{\alpha}}  y_2^{\vec {\alpha}} r_2^{i_2} z_2^{\vec {\alpha}} 
  \cdots y_k^{\vec {\alpha}} r_k^{i_k} z_k^{\vec {\alpha}}))\\
  &=\tau (e_2 \Phi (\sum_{\vec {\alpha} \in \{ 0 \}^k} y_1^{\vec {\alpha}} 
  z_1^{\vec{\alpha}}  y_2^{\vec {\alpha}} z_2^{\vec {\alpha}} 
  \cdots y_k^{\vec {\alpha}} z_k^{\vec {\alpha}}))\\
 &= \tau (e_2 \Phi (\sum_{\vec {\alpha} \in \{ 0 \}^k} \frac {1}{d^2} 
\sum_{m_3, m_4=0}^{d-1} 
U_{L_2}^{m_4} 
  U_{L_1}^{m_3}   (y_1^{\vec {\alpha}} z_1^{\vec{\alpha}}) {U_{L_1}^*}^{m_3} 
   {U_{L_2}^*}^{m_4}  y_2^{\vec {\alpha}} z_2^{\vec {\alpha}} 
  \cdots y_k^{\vec {\alpha}} z_k^{\vec {\alpha}}))\\
   &= \tau (e_2 \Phi (\sum_{\vec {\alpha} \in \{ 0 \}^k} x_1^{\vec {\alpha}} 
   y_2^{\vec {\alpha}} 
   z_2^{\vec {\alpha}} \cdots y_k^{\vec {\alpha}} z_k^{\vec {\alpha}}))=\\
  &\cdots\\
  &=
  \tau (e_2 \sum_{\vec {\alpha} \in \{ 0 \}^k} \Phi(x_1^{\vec {\alpha}})
  \cdot \Phi(x_2^{\vec {\alpha}}) \cdots \Phi(x_k^{\vec {\alpha}}))\\
   &=\lambda
  \tau (\sum_{\vec {\alpha} \in \{ 0 \}^k} \Phi(x_1^{\vec {\alpha}})
  \cdot \Phi(x_2^{\vec {\alpha}}) \cdots \Phi(x_k^{\vec {\alpha}}))\\
\end{align*}

\end{proof}

\begin{remark}
The trace $\tau$ on ${\mathfrak M}_m$ gives a hyperfinite Markov trace 
on the universal Jones algebra associated to $A_{1, m} \subset A_{0, m}$
\cite {sP}.
\end{remark}

\subsection {The $m=3$ Case}
\begin{align*}
   &A_{1, 3}= \{ e_3 \}''\\
   &={\mathbb C} e_3 \oplus {\mathbb C} (1- e_3)\\
   &\subset {\mathbb C} \otimes M_p ({\mathbb C})
   \oplus  {\mathbb C} \otimes M_{q-p} ({\mathbb C})\\
   &\subset M_q ({\mathbb C}) 
\end{align*}

There is an intermediate subalgebra,
\begin{align*}
   &A_{0, 3}= \{ e_2, e_3 \}''\\
   &= M_2 ({\mathbb C}) \oplus {\mathbb C}\\
   &\subset M_2 ({\mathbb C}) \otimes M_p ({\mathbb C})
   \oplus  {\mathbb C} \otimes M_{q-2p} ({\mathbb C})\\
   &\subset M_q ({\mathbb C}) 
\end{align*}

By the above theorem, there is an endomorphism of $\Phi$ of 
the hyperfinite $II_1$ factor such that:
\[
  \Phi^k (R)' \cap R= \otimes^k A_{1, 3}
\]

Take $Q= \Phi (R)$,
we have a tower of inclusions of hyperfinite $II_1$ factors
with the trace $\tau$:
\[
  Q \subset <Q, e_2> \subset <Q, e_2, e_3>
\]

The problem is to determine whether $<Q, e_2>$ is equal to
$<Q, e_2, e_3>$ or not.

Note that $[Q, e_3]= 0$. 
Let $\theta$ be the isomorphism of $Q e_3$ onto $Q (1- e_3)$.
Then any element in $Q$ can be decomposed as:
\[
  x_1 \oplus \theta (x_1), \quad x_1 \in Q e_3
\]

In the matrix form, we can write
\[
  e_3= \begin{bmatrix}1&0&0\\0&0&0\\0&0&0\end{bmatrix}, \quad
  e_2= \begin{bmatrix}\lambda&\sqrt{\lambda (1-\lambda)}&0\\
  \sqrt{\lambda (1- \lambda)}&1- \lambda&0\\0&0&0\end{bmatrix}
\]

For an element $x \in Q$, we can write
\begin{align*}
  &x= \begin{bmatrix}x_1&0&0\\0&x_2&x_{23}\\0&x_{32}&x_3\end{bmatrix}\\
  &\theta (x_1)= \begin{bmatrix}x_2&x_{23}\\x_{32}&x_3\end{bmatrix}\\
  &x_1= e_3 \cdot x \cdot e_3\\ 
  &x_2 = (e_2 \vee e_3- e_3) \cdot x \cdot (e_2 \vee e_3- e_3)\\
  &x_{23} = (e_2 \vee e_3- e_3) \cdot x \cdot (1- e_2 \vee e_3)\\
  &x_{32} = (1- e_2 \vee e_3) \cdot x \cdot (e_2 \vee e_3- e_3)\\
  &x_3 = (1- e_2 \vee e_3) \cdot x \cdot (1- e_2 \vee e_3)
\end{align*}

The matrix calculation shows:
\begin{align*}
   &e_2 x e_2=\\
   &e_2 \cdot \begin{bmatrix}\lambda x_1 + (1- \lambda) x_2&0&0\\
  0&\lambda x_1 + (1- \lambda) x_2&0\\0&0&x_3\end{bmatrix}\\ 
  &e_3 e_2 x e_2 x' e_2 e_3=\\
   &\lambda e_3 \cdot
   \begin{bmatrix}\lambda x_1 + (1- \lambda) x_2&0&0\\
  0&0&0\\0&0&0\end{bmatrix}
  \cdot \begin{bmatrix}\lambda x'_1 + (1- \lambda) x'_2&0&0\\
  0&0&0\\0&0&0\end{bmatrix}\\
\end{align*}

Note that
\begin{align*}
  &e_3 e_2 x e_2 e_3= \lambda e_3 \cdot 
  \begin{bmatrix}\lambda x_1 + (1- \lambda) x_2&0&0\\
  0&0&0\\0&0&0\end{bmatrix}\\
  &= \lambda e_3 \cdot 
  [(\lambda x_1 + (1- \lambda) x_2) \oplus
  \theta (\lambda x_1 + (1- \lambda) x_2)]
  \quad \in e_3 Q
\end{align*}

To summarize, $e_3$ implements a conditional expectation of
$<Q, e_2>$ onto $Q$.

We collect some useful facts. 
\begin{align*}
   &\forall x \in Q, \quad [x, e_3]= 0\\
   &\tau (e_3 x)= \tau ({\rm E}^{\tau}_Q (e_3) x) =\lambda \tau (x)\\
   &\tau (x_1)= \lambda \tau (x) \\
    &=\lambda \tau (x_1)+ 
   \lambda \tau (x_2)+ \lambda \tau (x_3)\\
   &\tau (e_2 x)= \lambda \tau (x)\\
   &\lambda \tau (x_1)+
    (1-\lambda) \tau (x_2) = \tau (x_1)\\
   &\tau (x_1)= \tau (x_2) \\
   &(1- 2 \lambda) \tau (x_1)=
   \lambda \tau (x_3) 
\end{align*}

By the above,
we have the  important identity:
\[
  {\rm E}^{\tau}_{<Q, e_2>} (e_3)= \lambda
\]
We also proved that the conditional expectation implemented 
by $e_3$ is trace-preserving.

We establish that
\[
  Q \subset <Q, e_2> \subset <Q, e_2, e_3>
\]
is a Jones basic construction with:
\begin{align*}
  &[<Q, e_2, e_3>: <Q, e_2>]= \lambda^{-1} \\
  &<Q, e_2>' \cap <Q, e_2, e_3>= {\mathbb C}
\end{align*}

\subsection {The $m=4$ Case}
\begin{align*}
   &A_{1, 4}= \{ e_3, e_4 \}''\\
   &=M_2 ({\mathbb C})  \oplus {\mathbb C} \\
   &\subset M_2({\mathbb C}) \otimes M_{pq} ({\mathbb C})
   \oplus  {\mathbb C} \otimes M_{q^2-2pq} ({\mathbb C})\\
   &\subset M_{q^2} ({\mathbb C}) 
\end{align*}

There is an intermediate subalgebra,
\begin{align*}
   &A_{0, 3}= \{ e_2, e_3, e_4 \}''\\
   &= M_2 ({\mathbb C}) \oplus M_3 ({\mathbb C}) \oplus {\mathbb C}\\
   &\subset M_2 ({\mathbb C}) \otimes M_{p^2} ({\mathbb C})
    \oplus M_3 ({\mathbb C}) \otimes M_{pq- p^2} ({\mathbb C})
   \oplus  {\mathbb C} \otimes M_{q^2-3pq+ p^2} ({\mathbb C})\\
   &\subset M_{q^2} ({\mathbb C}) 
\end{align*}

By the main theorem, there is an endomorphism of $\Phi$ of 
the hyperfinite $II_1$ factor such that:
\[
  \Phi^k (R)' \cap R= \otimes^k A_{1, 4}
\]

Take $Q= \Phi (R)$,
we have a tower of inclusions of hyperfinite $II_1$ factors
with the trace $\tau$:
\[
  Q \subset <Q, e_2> \subset <Q, e_2, e_3> \subset <Q, e_2, e_3, e_4>
\]

The problem is to determine whether $<Q, e_2, e_3, e_4>$ properly contains
$<Q, e_2, e_3>$. Nor is known whether $<Q, e_2>$ is a factor or not.

In the matrix form, we can write
\begin{align*}
  &e_4= \begin{bmatrix}1&0&0&0&0&0\\
       0&0&0&0&0&0\\0&0&1&0&0&0\\
       0&0&0&0&0&0\\0&0&0&0&0&0\\
       0&0&0&0&0&0\end{bmatrix}\\
  &e_3= \begin{bmatrix}\lambda&{\sqrt{\lambda (1- \lambda)}}&0&0&0&0\\
       {\sqrt{\lambda (1- \lambda)}}&1- \lambda&0&0&0&0\\
      0&0&\lambda&{\sqrt{\lambda (1- \lambda)}}&0&0\\
       0&0&{\sqrt{\lambda (1- \lambda)}}&1- \lambda&0&0\\0&0&0&0&0&0\\
       0&0&0&0&0&0\end{bmatrix}\\ 
  &e_2= \begin{bmatrix}1&0&0&0&0&0\\
       0&0&0&0&0&0\\0&0&0&0&0&0\\
       0&0&0&\frac{\lambda}{1- \lambda}&\frac{\sqrt{\lambda (1- 2 \lambda)}}
      {1- \lambda}&0\\
       0&0&0&\frac{\sqrt{\lambda (1- 2 \lambda)}}
      {1- \lambda}&\frac{1- 2 \lambda}{1- \lambda}&0\\
       0&0&0&0&0&0\end{bmatrix}
\end{align*}

For an element $x \in Q$, we can write
\[
  x= \begin{bmatrix} x_1&0&x_{12}&0&0&0\\
  0&x_1&0&x_{12}&0&0\\
  x_{21}&0&x_{2}&0&0&0\\
   0&x_{21}&0&x_{2}&0&0\\
   0&0&0&0&x_3&x_{34}\\
    0&0&0&0&x_{43}&x_4\\ 
  \end{bmatrix}
\]

The matrix calculation shows:
\begin{align*}
   &e_2 x e_2=\\
   &e_2 \cdot \begin{bmatrix}x_1&0&0&0&0&0\\
   0&0&0&0&0&0\\0&0&0&0&0&0\\
   0&0&0&\frac{\lambda}{1- \lambda} x_2+ \frac{1- 2\lambda}{1- \lambda} 
  x_3&0&0\\
   0&0&0&0&\frac{\lambda}{1- \lambda} x_2+ \frac{1- 2\lambda}{1- \lambda} 
  x_3&0\\0&0&0&0&0&0\end{bmatrix} 
\end{align*}

Take $({\mathfrak q}_1)_4$ as an example.
\begin{align*}
   &({\mathfrak q}_1)_4=\\
   &\begin{bmatrix} x_1&0&0&0&0&0\\
  0&x_1&0&0&0&0\\
  0&0&x_{2}&0&0&0\\
   0&0&0&x_{2}&0&0\\
   0&0&0&0&x_3&0\\
    0&0&0&0&0&x_4\\ 
  \end{bmatrix} \otimes {\bf 1} \otimes {\bf 1} \otimes 
  {\mathfrak q}_1\\
  &x_1= [1 \ \rho_1^{1} \ \rho_1^{2} \ \cdots \ \rho_1^{p^2-1}] 
   \in M_{p^2} ({\mathbb C})\\
  &x_2= [\rho_1^{p^2} \ \rho_1^{p^2+1} \ \rho_1^{p^2+2} \ \cdots \ \rho_1^{pq-1}]
  \in M_{pq-p^2} ({\mathbb C})\\
  &x_3= {\bf 1}_{\in M_{pq-p^2} ({\mathbb C})}\\
  &x_4= {\bf 1}_{\in M_{q^2- 3pq +p^2} ({\mathbb C})}
\end{align*}

Note that
\begin{align*}
   &{\rm E}^{\tau}_{A_{1, 4}} (e_2 ({\mathfrak q}_1)_4 e_2 
    {({\mathfrak q}_1)_4}^* e_2) \\
   &= 2 \lambda \begin{bmatrix}\lambda&0&0\\0&\kappa&0\\0&0&\kappa\end{bmatrix}\\
   &= 2 \lambda^2 e_4+ 2 \lambda \kappa (1- e_4)\\
   &{\rm where \ } \kappa < \lambda (1- \lambda)
\end{align*}

By the relative Dixmier property \cite{sP99}:
\[
  {\rm E}^{\tau}_{A_{1, 4}} (e_2 ({\mathfrak q}_1)_4 e_2 
  {({\mathfrak q}_1)_4}^* e_2)=
  {\rm E}^{\tau}_{Q' \cap R}  (e_2 ({\mathfrak q}_1)_4 e_2 
  {({\mathfrak q}_1)_4}^* e_2)
  \in <e_2, Q>
\]

Therefore 
\begin{align*}
    &e_4 \in <Q, e_2>\\
    &Z(<Q, e_2>)= {\mathbb C} e_4 \oplus {\mathbb C} (1- e_4)\\
    &<Q, e_2, e_3>= <Q, e_2, e_3, e_4>
\end{align*}

However we can consider the inclusion of hyperfinite $II_1$ factors:
\[
  Q e_4 \subset <Q, e_2> e_4
\]

An element in $Q e_4$ can be written as
\[
  \begin{bmatrix}x_1&x_{12}\\x_{21}&x_2\end{bmatrix}
\]

Since 
\[
  e_2 e_4= \begin{bmatrix}1&0\\0&0\end{bmatrix},
\]
$<Q, e_2> e_4$ is generated by:
\[
  \begin{bmatrix}x_1&0\\0&0\end{bmatrix}, \ 
  \begin{bmatrix}0&x_{12}\\0&0\end{bmatrix}, \ 
   \begin{bmatrix}0&0\\x_{21}&0\end{bmatrix}, \
   \begin{bmatrix}0&0\\0&x_2\end{bmatrix}
\]

$e_2 e_4$ implements a conditional expectation from $Q e_4$ onto
$\{ e_2 e_4 \}' \cap Q e_4$ by
\begin{align*}
  &\begin{bmatrix}1&0\\0&0\end{bmatrix} \cdot
  \begin{bmatrix}x_1&x_{12}\\x_{21}&x_2\end{bmatrix} 
   \cdot \begin{bmatrix}1&0\\0&0\end{bmatrix}\\
   &=\begin{bmatrix}1&0\\0&0\end{bmatrix} \cdot 
   \begin{bmatrix}x_1&0\\0&x_2\end{bmatrix}
\end{align*}
Note that 
\[ 
  \begin{bmatrix}x_1&x_{12}\\x_{21}&x_2\end{bmatrix} \in Q e_4
  \Rightarrow
  \begin{bmatrix}x_1&0\\0&x_2\end{bmatrix} \in Q e_4
\]
because of the Dixmier property.


By the above,
we have the  important identity:
\[
  {\rm E}^{\tau}_{Q e_4} (e_2 e_4)= \lambda e_4
\]
We also proved that the conditional expectation implemented 
by $e_2 e_4$ is trace-preserving.

We establish that
\[
  \{ e_2 e_4 \}' \cap Q e_4 \subset Q e_4 \subset <Q, e_2> e_4
\]
is a Jones basic construction with:
\begin{align*}
  &[<Q, e_2> e_4: Q e_4]= \lambda^{-1} \\
  &Q e_4' \cap <Q, e_2> e_4 = {\mathbb C} e_4
\end{align*}

\section {Analogy} \label {S:anal}
In this section, we construct inclusions of non-hyperfinite
$II_1$-factors via free product with amalgamation as a comparison
to the main theorem.

\begin {theorem}
For any finite dimensional $C^*$-algebra $A$
with any trace vector
$\vec s$ whose components are rational numbers, there exists a tower
of inclusions of $II_1$-factors, $M \subset M_1 \subset M_2 \subset
M_3 \subset \cdots$, with the trace $\tau$ such that
\[
   M' \cap M_k = \otimes_{i=1}^k A.
\]
The canonical trace $\tau$ on $M_k$ extends the trace vector $\vec s$
on $A$. 
\end {theorem}
The main tool is the relative commutant theorem by S.Popa
\cite {sP}.

\begin {lemma} \cite {sP}
Let $(P_1, \tau_1), (P_2, \tau_2)$ be two finite von Neumann
algebras with a common von Neumann subalgebra $B \subset P_1$, $B
\subset P_2$, such that $P_1= Q {\overline {\otimes}} B$ where $Q$
is a nonatomic finite von Neumann algebra. If $(P, \tau)$ denotes
the amalgamated free product $(P_1, \tau_1)*_B (P_2, \tau_1)$ then
${Q_0}' \cap P = ({Q_0}' \cap Q) \overline {\otimes} B$ for any
nonatomic von Neumann subalgebra $Q_0 \subset Q$.

Assume that $L \subset P$ is a von Neumann subalgebra satisfying
the properties:\\
(1) $Q_0 \subset L$.\\
(2) $L \cap P_2$ contains an element $y \not= 0$ orthogonal to $B$,
i.e. $E_B (y)= 0$ and with $E_B (y^* y) \in \mathbb {C} 1$.\\
Then $L' \cap P= L' \cap B$. If in addition $L \cap B= \mathbb {C}$
then $L$ is a type $II_1$ factor.
\end {lemma}

We can embed $A$ in the full matrix algebra $M_d (\mathbb
{C})$. Consider the tensor product of $M_d (\mathbb {C}) \otimes M_2
(\mathbb {C})$. Identify $x \in A$ as $x \otimes {\bf 1}_{M_2
(\mathbb {C})}.$
Take the element $y$:
\[
  y= {\bf 1}_{M_d (\mathbb {C})} \otimes {\begin {pmatrix}
       0 & 1\\
       1 & 0
      \end {pmatrix}}.
\]
In $M_{d} (\mathbb {C}) \otimes M_{2} (\mathbb {C})$, 
$y \not= 0$ is an element orthogonal to
$A$, i.e., $E_{A} (y)= 0$ and with $E_{A} (y^* y)= 1$.

Let $M$ be a  $II_1$-factor. We construct $M_1$ via the
following map $\Gamma$:
\[
  M_1 = \Gamma (M) = (M \otimes A)*_A M_{2n} (\mathbb {C}).
\]
The trace $\tau$ on $M$ can be extended to $M_1$.

Put
\begin {align*}
  &P_1 = M \otimes A\\
  &P_2 = M_{2d} (\mathbb {C})= M_{d} (\mathbb {C}) 
  \otimes M_{2} (\mathbb {C})\\
  &L= P= M_1
\end {align*}

We get:
\[  
  {M_1}' \cap M_1= {M_1}' \cap A \subseteq {M_{2d} (\mathbb {C})}'
  \cap A = \mathbb {C}
\]
That is, $M_1$ is a nonhyperfinite $II_1$ factor.

The relative commutant
\[
  M' \cap M_1= (M' \cap M) \otimes A= A
\]

Viewing $\Gamma$ as a machine producing nonhyperfinite $II_1$ factors, we get an
ascending towers of $II_1$ factors:
\[
  M \subset \Gamma (M)= M_1 \subset \cdots M_i
  \subset \Gamma (M_i)= M_{i+1} \cdots.
\]
There is a unique trace $\tau$ on every $M_i$.

We calculate the relative commutant $M' \cap M_k$ by
induction.
Assume
\[
  M' \cap M_k= \otimes^k A.
\]
By the above lemma,
\begin {align*}
  &M' \cap M_{k+1}= (M' \cap M_k) \otimes A\\
  &= (\otimes^k A) \otimes A= \otimes^{k+1} A
\end {align*}

In the end, we would boldly suggest an analogy between binary shifts
and free product with amalgamation.

\begin {acknowledgments}
I wish to express my gratitude toward V.Jones for proposing the problem
and for helpful advice. I am indebted to the following persons for correcting 
my mistakes and for useful discussions: M.Choda, M.Izumi, Y.Kawahigashi, Z.Landau,
S.Neshveyev, G.Powers, R.Price, E.St{\o}rmer. I would like to thank Jing Yu for his 
hospitality during my stay at the National Center for Theoretical Sciences in Taiwan.
\end {acknowledgments}


\begin {thebibliography} {9}
\bibitem {dB}
D.Bisch, \emph{Bimodules, higher relative commutants and the fusion
algebra associated to a subfactor}, Operator algebras and their
applications (Waterloo, ON, 1994/1995), 13-63, Fields Inst. Commun.,
13, Amer. Math. Soc., Providence, RI, 1997
\bibitem {mC87}
M.Choda, \emph {Shifts on the hyperfinite $II_1$ factor},
J. Operator Theory, 17 (1987),223-235
\bibitem{vJ}
V.F.R.Jones, \emph{Index for subfactors},  Invent. Math.  72  (1983), no. 1, 1-25
\bibitem {hN95}
H.Narnhofer, W. Thirring, E.St{\o}rmer, \emph{$C^*$-dynamical
systems for which the tensor product formula for entropy fails},
Ergodic Theory Dynam. Systems 15 (1995), no. 5, 961-968
\bibitem {rP88}
R.T.Powers, \emph {An index theory for semigroups of
$*$-endomorphisms of $\mathbb {B} (\mathcal {H})$ and type $II_1$
factors}, Canad. J. Math, 40 (1988), 86-114
\bibitem {sP83}
S.Popa,  \emph{Orthogonal pairs of $*$-subalgebras in finite von
Neumann algebras}, J. Operator Theory 9 (1983), no. 2, 253-268
\bibitem {sP}
S.Popa, \emph {Markov traces on universal Jones algebras
and subfactors of finite index}, Invent. Math, 111(1993),
no 2, 375-405
\bibitem {sP94}
S.Popa, \emph {Classification of amenable subfactor of
type $II$}, Acta. Math, 172(1994),
no 2, 163-255
\bibitem {sP95}
S.Popa, \emph{An axiomatization of the lattice of higher relative
commutants of a subfactor}, Invent. Math. 120 (1995), no. 3, 427-445
\bibitem {sP99}
S.Popa, \emph{The relative Dixmier property for inclusions of von Neumann
 algebras of finite index},  Ann. Sci. Ecole Norm. Sup. (4)  32  (1999),  
no. 6, 743-767
\end {thebibliography}
\end {document}